\documentclass[12pt]{article}
\usepackage{pictexwd,dcpic} 
\usepackage[small,nohug,heads=littlevee]{diagrams}
\diagramstyle[labelstyle=\scriptstyle]
\usepackage{amsmath}
\usepackage{amsfonts}
\usepackage{epsfig}
\usepackage[dvips]{color}
\usepackage{epsfig}
\usepackage{theorem}
\newtheorem{theorem}{Theorem}[section]

\newtheorem{lemma}{Lemma}[section]

\newtheorem{definition}{Definition}[section]

\newtheorem{proposition}{Proposition}[section]

\newtheorem{corollary}{Corollary}[section]

\def \proof {\vspace{-.4cm}{\bf Proof \quad}}

\def \C {\mathbb C}

\def \cbar {{\overline {\mathbb C}}}

\def \D {\mathbb D}

\def \DO {({\mathbb D},0)}

\def \N {\mathbb N}

\def \T {\mathbb T}

\def \Car {Carath\'eodory }

\def \diams {{\rm diam}^\#}

\def \Do {({\mathbb D}, 0)}

\def \epsilon {\varepsilon}

\def \phi {\varphi}

\def \pt {{\rm pt}}

\def \Am {{\cal A}_{\infty, m}}

\def \Fm {{\cal F}_m}

\def \Fn {{\cal F}_n}

\def \J {{\cal J}}

\def \Jm {{\cal J}_m}

\def \Jmn {{\cal J}_m^n}

\def \Jn {{\cal J}_n}

\def \K {{\cal K}}

\def \Km {{\cal K}_m}

\def \Kmp {{\cal K'}_m}

\def \Pm {\{P_m \}_{m=1}^\infty}

\def \Pmn {\{P_m^n \}_{m=1}^\infty}

\def \Aa {(A,a)}

\def \Aaaa {(A_\alpha, a_\alpha)}

\def \Aaaap {(A'_\alpha, a'_\alpha)}

\def \bfcF {\boldsymbol {\check {\mathcal F}}}

\def \bfg {{\boldsymbol \gamma}}

\def \bfG {{\boldsymbol \Gamma}}

\def \bfcg {{\boldsymbol {\check \gamma}}}

\def \bfcG {{\boldsymbol {\check \Gamma}}}

\def \bftg {{\boldsymbol {\tilde \gamma}}}

\def \bftG {{\boldsymbol {\tilde \Gamma}}}

\def \clm {\check{l}_m}

\def \cLm {\check{L}_m}

\def \tLm {{\tilde L}_m}

\def \gm {\gamma_m}

\def \cgm {\check{\gamma}_m}

\def \tgm {{\tilde \gamma}_m}

\def \ga {\gamma_\alpha} %%%%%From Annuli2 %%%%%%%%

\def \Gm {\Gamma_m}

\def \cGm {\check{\Gamma}_m}

\def \tGm {{\tilde \Gamma}_m}

\def \tRm {\tilde R_m}

\def \Uu {(U,u)}

\def \Uaua {(U_\alpha, u_\alpha)}

\def \Uauap {(U'_\alpha, u'_\alpha)}

\def \Umum {(U_m,u_m)}

\def \Unun {(U_n,u_n)}

\def \Vava {(V_\alpha, v_\alpha)}

\def \Wawa {(W_\alpha, w_\alpha)}

\def \cpsi {\check{\psi}}

\def \ccpsi {{\check{\check \psi}}}

\def  \fuv {(f: U \rightarrow V) }

\def  \fuvm {\{f_{m+1}: U_m \rightarrow V_{m+1}\}_{m=0}^\infty}

\def  \fuvmp{\{f_{m+1}: U'_m \rightarrow V'_{m+1}\}_{m=0}^\infty}

\def  \fuvmt {\{\tilde f_{m+1}: \tilde U_m \rightarrow \tilde V_{m+1}\}_{m=0}^\infty}

\def \calA {{\mathcal A}}

\def \calAp {{\mathcal A}'}

\def \calF {{\mathcal F}}

\def \calG {{\mathcal G}}

\def \calL {{\mathcal L}}

\def \calLc {\check{\mathcal L}}

\def \calLt {\tilde {\mathcal L}}

\def \calM {{\mathcal M}}

\def \calS {\mathcal S}

\def \calT {\mathcal T}

\def \calU {\mathcal U}

\def \calUp {{\mathcal U}'}

\def \calV {\mathcal V}

\def \calVp {{\mathcal V}'}

\def \calW {\mathcal W}

\def \Psic {\check{\Psi}}

\def \Psicc {\check{\check{\Psi}}}

\def \bfccPsi {{\boldsymbol {\check {\check \Psi}}}}

\def \FUV {(\calF : \calU \rightarrow \calV)}

\def \FUVa {\{(\calF^{\alpha} : \calU^{\alpha} \rightarrow \calV^{\alpha})\}_{\alpha \in A}}

\def \FUVn {\{(\calF^n : \calU^n \rightarrow \calV^n)\}_{n \ge 1}}

\def \FUVp {(\calF : \calU' \rightarrow \calV')}

\def \FUVt {(\tilde \calF: \tilde \calU \rightarrow \tilde \calV)}

\title{A Straightening Theorem for non-Autonomous Iteration}

\author{Mark Comerford\\ 
Department of Mathematics\\
University of Rhode Island\\
5 Lippitt Road, Room 102F\\
Kingston, RI 02881, USA\\
email: {\tt mcomerford@math.uri.edu}}

\begin{document}  

\maketitle

\begin{abstract}The classical straightening theorem as proved by Douady and Hubbard states that a polynomial-like mapping is hybrid equivalent to a polynomial. 
We generalize this result to non-autonomous iteration where one considers composition sequences arising from a sequence of functions which in general is allowed to vary. In order to do this, new techniques are required to control the distortion and quasiconformal dilatation of the hybrid equivalence. In particular, the \Car topology for pointed domains allows us to specify the appropriate bounds on the sequence of sets on which the polynomial-like mapping sequence is defined and give us good estimates on the degree of distortion and quasiconformality.
\end{abstract}

\section{Basic Definitions and Results}

In their famous paper `On the Dynamics of Polynomial-like Mappings', Douady and Hubbard \cite{DH} introduced the concept of a polynomial-like mapping and used it to show the existence of small copies of the Mandelbrot set. One of the first steps of their proof was the straightening theorem where they showed a polynomial-like mapping is hybrid equivalent to a polynomial. 

In this paper, we extend their result to the setting of non-autonomous polynomial iteration where one composes a sequence of mappings which in general are allowed to vary. We generalize the concept of a polynomial-like mapping to a polynomial-like mapping sequence and show that such a sequence is hybrid equivalent in the appropriate sense to a sequence of polynomials which satisfies suitable bounds on the degrees and coefficients. This requires some new techniques, especially the use of the \Car topology for pointed domains and the notion of a bounded family of pointed domains where connectivity and conformal properties are preserved under limits in this topology.

\subsection{Bounded Polynomial Sequences}

Let $d \ge 2$, $M \ge 0$, $K \ge 1$ and let $\Pm$ be a sequence of polynomials where each $P_m(z) = a_{d_m,m}z^{d_m} + a_{d_m-1,m}z^{d_m-1} + \cdots \cdots + 
a_{1,m}z + a_{0,m}$ is a polynomial of degree $2 \le d_m \le d$ whose coefficients satisfy
\[\qquad 1/K \le |a_{d_m,m}| < K,\: m \ge 1, \quad |a_{k,m}| \le M,\: m \ge 1,\: 0 \le k \le d_m -1. \]
Such sequences are called \emph{bounded sequences of polynomials} or simply \emph{bounded sequences}. 

For each $1 \le m$, let $Q_m$ be the composition $P_m \circ \cdots \cdots \circ P_2 \circ P_1$ and for each $0 \le m < n$, let $Q_{m,n}$ be the composition $P_n \circ \cdots \cdots \circ P_{m+2} \circ P_{m+1}$. Let the degrees of these compositions be $D_m$ and $D_{m,n}$ respectively so that $D_m = \prod_{i=1}^m d_i$, $D_{m,n} = \prod_{i=m+1}^n d_i$.

For each $m \ge 0$ define the \emph{$m$th iterated Fatou set} or simply the \emph{Fatou set} at time $m$, $\Fm$, by 
\[ \Fm = \{z \in \cbar : \{Q_{m,n}\}_{n=m}^\infty \;
\mbox{is a normal family on some neighbourhood of}\; z \}\]
where we take our neighbourhoods with respect to the spherical topology on $\cbar$ 
and let the \emph{$m$th iterated Julia set} or simply the \emph{Julia set} at time $m$, $\Jm$, to be the complement $\cbar \setminus \Fm$. 

It is easy to show that these iterated Fatou and Julia sets are completely invariant in the following sense. 

\begin{theorem}For any $m \le n \in \N$, $Q_{m,n}(\Jm) = \Jn$  
and $Q_{m,n}(\Fm) = \Fn$, with Fatou 
components of $\Fm$ being mapped surjectively onto those of $\Fn$ by 
$Q_{m,n}$. 
\end{theorem}

If $\Pm$ is a bounded sequence, we can clearly find some radius $R$ 
depending only on the bounds $d$, $K$, $M$ above so that for any sequence $\Pm$ as above and any $m \ge 0$, it is easy to see that 
\[ |Q_{m,n}(z)| \to \infty \qquad \mbox{as} \quad n \to \infty, \quad |z| > R\]
which shows in particular that as for classical polynomial Julia sets, there 
will be a \emph{basin at infinity at time $m$}, $\Am$ on which all points escape to infinity under
iteration. Such a radius will be called an {\it escape radius} for the bounds $d$, $K$, $M$. Note that the maximum principle shows that just as in 
the classical case, there can be only one component on which $\infty$ is 
a normal limit function and so the sets $\Am$ are 
completely invariant in the sense given in Theorem 1.1.

The complement of $\Am$ is called the {\it filled Julia set} at time $m$ for
the sequence $\Pm$ and is denoted by $\Km$. As above, the same 
argument using Montel's theorem as in the classical case shows that 
$\partial \Km = \Jm$ as one might expect.

%%%%%%%%%%%%%%%%%%%%%%%%%%%%%%%%%%%%%%%%%%%%%

\section{Some Useful Lemmas}

In this section, we collect some results, mostly concerning distortion of anlaytic mappings, which 
will be of use to us later on. Firstly, in our proof of the straightening theorem for non-autonomous iteration, we 
will be making extensive use of univalent mappings. Hence the Koebe one-quarter theorem and 
the other standard distortion results for univalent functions will naturally be of great importance to 
us. However, as this material is standard and well known, we will merely point the reader in the
direction of some of the many excellent references on the subject (e.g. \cite{Dur, CG, Po2}).
Less well-known but also extremely useful is the result below. Roughly, what it states is that for 
mappings of bounded degree, the sizes of discs in the hyperbolic metric cannot shrink by too much. 
A proof can be found in \cite{CJY} where the authors make extensive use of the result in their 
investigation of semi-hyperbolic behaviour. 
 
\begin{lemma} Let $D \subset \C$ be a simply connected domain and let
$F :  D \to \D$, $F(\partial D) \subset \partial \D$, be $p$-valent 
(i.e. of topological degree $p$). Then if $\rho$ denotes the hyperbolic metric, 
there exists a positive constant $C$ which depends only on $p$ such that 
\begin{eqnarray*}
\{w \in \D: \rho_\D(F(z_0),w) \le 1/C\} &\subset& 
F(\{z \in D: \rho_D(z,z_0) \le 1\})\\
&\subset& \{w \in \D: \rho_\D(F(z_0),w) \le 1\}.
\end{eqnarray*}
\end{lemma}

Next we recall (e.g. \cite{CG} page 15 Theorem 4.3) that for a simply connected domain $U$, we have the following estimate for the 
hyperbolic metric:

\[ \frac{|dz|}{2\delta(z)} \le d\rho(z) \le \frac{2|dz|}{\delta(z)},\]

where $\delta(z)$ denotes the Euclidean distance from $z$ to the boundary of $U$.
Here, the upper bound is valid for all domains while the lower bound only holds in the  simply 
connected case. Nonetheless, one can completely characterize all domains for which one 
has a lower bound of the same form as that given above. Recall that a set $X$ is $K$-uniformly perfect
if the moduli of annuli which separate the set (i.e. each of the complememtary components of the 
annulus contains points of $X$) are bounded above by $K$. One can show (e.g. \cite{Po1}) that those domains for which one has a lower estimate of the same form as above are precisely those whose boundaries are uniformly perfect. An obvious special case is a conformal annulus. 

\begin{corollary} 
Let $A$ be an annulus of modulus $M$. Then there exists\\ $c=c(M) > 0$ (which is decreasing in $M$) such that if $\rho(.\,,.)$ denotes the 
hyperbolic metric on $A$, then 

\[ \frac{c|dz|}{\delta(z)} \le d\rho(z) \le \frac{2|dz|}{\delta(z)}.\]

\end{corollary}

%%%%%%%%%%%%%%%%%%%%%%%%%%%%%%%%%%%%%%%%%%%%%

\section{The \Car Topology}

In non-autonomous iteration, one deals constantly with questions of convergence for 
varying sequences of sets and varying sequences of functions defined on these 
sets. The \Car topology is ideally suited to dealing with this 
situation and provides the perfect setting to make the results in the later
sections of this paper as natural and simple as possible to state and prove. 

The \Car topology is most commonly used for discs, but it also works for domains of arbitrary connectivity. 
A {\it pointed domain} is a pair $(U,u)$ consisting of an open  
connected subset $U$ of $\cbar$, (possibly equal to $\cbar$ itself) and a point $u$ in 
$U$. We say that $\Umum \to \Uu$ in the \Car topology if and only if

\begin{enumerate}

\item $u_m \to u$ in the spherical topology;

\item for all compact sets $K \subset U$, $K \subset U_m$ for all but finitely many $m$;

\item for any {\it connected} open set $N$ containing $u$, if $N \subset U_m$ for
infinitely many $m$, then $N \subset U$.

\end{enumerate}

We also wish to consider the degenerate case where $U = \{u\}$. In this 
case, condition 2. is omitted ($U$ has no interior of which we can take
compact subsets) while condition 3. becomes

\begin{description}

\item[{\quad \rm 3.}]  for any {\it connected} open set $N$ containing $u$, 
$N$ is contained in at most finitely many of the sets $U_m$.

\end{description}

The above definition is a slight modification of that given in the book of McMullen \cite{McM} and much 
of what follows in this section is based on his exposition.  It is not too hard to show the following useful result. 

\begin{theorem}
Let $\Unun$ be a sequence of pointed domains and let $\Uu$ be a pointed domain (where we include the degenerate case where $U = \{u\}$). Then $\Unun \to \Uu$ in the \Car topology if and only if $u_n \to u$ and for any subsequence where the complements of the sets $U_n$ 
converge in 
the Hausdorff topology (with respect to the spherical metric), $U$ correspsonds with the connected component of the
complement of the Hausdorff limit which contains $u$ (this component being 
empty in the degenerate case $U = \{u\}$).
\end{theorem}

\Car convergence can also in terms of the \emph{\Car kernel} (e.g. \cite{Dur}). Lastly, for hyperbolic domains (i.e. domains whose universal covering space is the unit disc), Epstein in his Ph.D. thesis \cite{Ep} showed that it is also equivalent to the local uniform convergence of suitably normalized covering maps on the unit disc. 

\Car convergence allows us to give a definition for convergence of holomorphic 
functions defined on sequences of pointed domains. This definition is given in \cite{Ep} and is a slight modification of that given by McMullen. Note that we use the notation ${\mathrm d^\#}(x,y)$ to denote the spherical distance between points $x$, $y$ in $\cbar$.

\begin{definition}
Let $\{\Unun\}_{n=1}^\infty$ be a sequence of pointed domains which converges in the \Car topology to a pointed domain $\Uu$ with $\Uu \ne (\{u\},u)$. For each $n$ let $f_n$ be an analytic function (with respect to the spherical topology) defined on $U_m$ and let $f$ be an analytic function defined on $U$. We say that $f_m$ converges to $f$ {\rm uniformly on compact subsets of $U$} or simply {\rm  locally uniformly to $f$ on $U$} if, for every compact subset $K$ of $U$ and every $\epsilon >0$, there exists $n_0$ such that ${\mathrm d^\#}(f_n(z), f(z)) < \epsilon$ on $K$ for all $n \ge n_0$.
\end{definition}

A fact worth noting is that the connectivity of the pointed domains cannot increase under \Car convergence. To be precise, if each $U_n$ is at most $N$-connected, then so is the \Car limit $U$. The reason for this is that complementary components are allowed to merge in the limit, but they cannot split up into more components. For example, the pointed conformal annuli $(\D \setminus \overline {\mathrm D}((1-\tfrac{2}{n}), \tfrac{1}{n}), 0)$, $n \ge 2$ are easily seen to converge to $(\D, 0)$.

%%%%%%%%%%%%%%%%%%%%%%%%%%%%%%%%%%%%%%%%%%%%%%%%%%%

\section{Equicontinuity and Bi-Equicontinuity}

In this section we continue collecting results dealing with sequences
of functions on variable domains which will be of use to us later. The proofs of all of these may be found in \cite{Com4} (for domains of arbitrary finite connectivity rather than just discs or annuli as we consider here). We 
begin with some definitions. It will be convenient for us to also consider degenerate families of domains where we either just have points 
or the whole Riemann sphere and these we will denote by $\pt$ and $\cbar$ respectively. A finitely connected domain $U \subset \cbar$ is called \emph{non-degenerate} if none of its complementary components is a point. 

\begin{definition} Let $n \ge 1$ and let $\calU = \{\Uaua\}_{\alpha \in A}$ be a family of 
pointed hyperbolic domains. 

We say that $\calU$ is {\rm hyperbolically non-degenerate} or simply {\rm non-degenerate} if every sequence in $\calU$ which is convergent in the \Car topology has a limit which is a pointed hyperbolic domain. 

We say that a family $\calU$ of $n$-connected non-degenerate pointed domains is {\rm bounded} if every sequence in $\calU$ which is convergent in the \Car topology has a limit which is a non-degenerate $n$-connected pointed domain and we write $\pt \sqsubset \calU \sqsubset \cbar$. Otherwise we say that $\calU$ is {\rm unbounded}.
\end{definition}

Using the Hausdorff version of \Car convergence, it is clear that any family of pointed domains (whose connectivities in general can differ) is precompact in the 
sense that any sequence will have a convergent subsequence. The content of the definition of boundedness, then, is that any limit is non-trivial and non-degenerate and that connectivity is preserved. Note that if ${\cal A} = {\Aaaa}_{\alpha \in A}$ is such a family of annuli, then $0 < \inf_\alpha \bmod A_\alpha \le \sup_\alpha \bmod A_\alpha < \infty$.

The boundedness of a family of pointed discs or annuli can be expressed in a quantitive manner using the \emph{\Car bound}. For a point $z$ in an open subset $U$ of $\cbar$, we let $\delta^\#(U)(z)$ denote the spherical distance to the boundary of $U$ while for some set $E \subset \cbar$, we let $\diams E$ denote the diameter of $E$ in the spherical metric.

\begin{definition} Let $\calA = \{\Aaaa\}_{\alpha \in A}$ be a family of pointed annuli. Let $\rho_{A_\alpha} (.,.)$ denote the hyperbolic metric for each annulus $A_\alpha$, let $\ga$ be the equator of this annulus and let $\ell_{A_\alpha}(\ga)$ denote its hyperbolic length in $A_\alpha$. 
We define the \emph{\Car bound} of $\calA$, $||\calA||$ by

\vspace{-.4cm}
\[ ||\calA|| = \sup_\alpha |\log ( \delta^\#_{A_\alpha}(a_\alpha)\, \diams (\cbar \setminus A_\alpha)) | \,+\, \sup_\alpha  \ell_{A_\alpha} (\ga) \,+\, \sup_\alpha \rho_{A_\alpha} (a_\alpha, \ga).\]

where $\rho_{A_\alpha} (a_\alpha, \ga)$ is the hyperbolic distance in $A_\alpha$ from $a_\alpha$ to $\ga$. For a family $\calU = \{\Uaua\}_{\alpha \in A}$ of pointed discs, the \emph{\Car bound} of $\calU$, $|| \,\calU ||$, is defined similarly where the second and third terms above are zero (there being no closed geodesics in this case).
\end{definition}

\vspace{-.2cm}
One can think of the three terms above as each preventing a
different way in which a \Car limit of non-degenerate annuli can fail to be another non-degenerate annulus. The first term prevents any \Car limit of pointed annuli (or discs) from simply being a point or $\cbar$ with just one point removed (recall that the spherical diameter of all of $\cbar$ is $\tfrac{\pi}{2}$ and in particular finite). The second term prevents complementary components from merging or becoming points (as in the example given at the end of the last section) while the third term prevents components of the complement from being `engulfed' by other components (\cite{Com4} gives an example of what can happen when these terms are allowed to become unbounded).  The following and Theorems 4.3 and 4.5 are special cases of parts of a larger result in \cite{Com4}.

\begin{theorem}[\cite{Com4} Theorem 4.2]
Let $\calU$ be a family of pointed discs or a family of non-degenerate pointed annuli. Then $\pt \sqsubset \calU \sqsubset \cbar$ if and only if $||\, \calU || < \infty$.
\end{theorem}

Since the equator of an annulus gives us a natural frame of reference, we will often specify pointed annuli without giving a base point in which case it is understood that the base point lies somewhere on the equator. Note that if we are dealing with a bounded family, it will not matter just where on the equator of each annulus we pick our base point. 

\begin{definition}
Let $\calA = \{\Aaaa\}_{\alpha \in A}$ be a family of pointed annuli (not necessarily bounded). We say that another family of 
pointed annuli, $\calAp = \{\Aaaap\}_{\alpha \in A}$ is {\rm bounded
above and below} or just {\rm bounded in $\calA$} if there 
exists $K \ge 1$ such that for every $\alpha \in A$, if we let $\gamma_\alpha'$ be the equator of each $A_\alpha'$, we have

\begin{enumerate}
\item $A'_\alpha$ is a subset of $A_\alpha$ which lies within hyperbolic distance at most $K$ about $a_\alpha$ in $A_\alpha$,

\item \[ \delta^\#_{A_\alpha'} (a'_\alpha) \ge \frac{1}{K}  \delta^\#_{A_\alpha} (a'_\alpha),\]

\item \[ \ell_{A_\alpha'}(\ga') \le K,\]

\item \[ \rho_{A_\alpha'} (a_\alpha', \ga') \le K.\]

\end{enumerate}

In this case, we write $\pt \sqsubset \calAp
\sqsubset \calA$.

Let $\calU = \{\Uaua\}_{\alpha \in A}$ be a family of pointed discs (not necessarily bounded in the 
sense given above). We say that another family of 
pointed discs, $\calUp = \{\Uauap\}_{\alpha \in A}$ is {\rm bounded
above and below} or just {\rm bounded in $\calU$} if there 
exists $K \ge 1$ such that for every $\alpha \in A$ the corresponding analogue of conditions 1. and 2. above hold (the remaining two conditions being vacuously true).
 \end{definition}

For the sake of brevity, from now on, unless otherwise specified, we will use 
$\calU = \{\Uaua\}_{\alpha \in A}, \calV = \{\Vava\}_{\alpha \in A}$ and $\calW = \{\Wawa\}_{\alpha \in A}$ to refer to bounded families of pointed discs or annuli where it is understood that all three families are indexed over the same set and all the domains of all three families possess the \emph{same} connectivity.

We are also interested in boundedness for the degenerate case of a family of simple closed curves, in our case, the equators of a family of conformal annuli. Recall that for a function $f$ defined on a subset of the Riemann sphere, we have the spherical derivative $f^\#(z) = \tfrac{f'(z)}{1 + |f(z)|^2}$. 

\begin{definition} Let $\Gamma = \{(\gamma_\alpha, z_\alpha)\}_{\alpha \in A}$ be a family of pointed simple smooth (i.e. ${\mathrm C}^1$) closed curves where for each $\alpha$, $z_\alpha$ is a point on $\gamma_\alpha$ and $\phi_\alpha : \T \mapsto \cbar$ is a parametrization of $\gamma_\alpha$. 

Let $\calU = \{\Uaua\}_{\alpha \in A}$ be a family of pointed domains of the same finite connectivity. We say $\Gamma$ is {\rm bounded above and below} 
or just {\rm bounded} in $\calU$ with constant $K \ge 1$ if 
\begin{enumerate}
\item For each $\alpha$, $\gamma_\alpha$ is a subset of $U_\alpha$ which lies within hyperbolic distance at most $K$ of $u_\alpha$ in $U_\alpha$; 

\vspace{.2cm}
\item The mappings $\phi_\alpha$ can be chosen so that the resulting family 
$\Phi$ is bi-equicontinuous on $\T$ in the sense that we can find $K \ge 1$ such that 

\[ \frac{1}{K} \le \frac{|\phi_\alpha^\#(z)|}{\delta_{U_\alpha}^\#(z_\alpha)} \le K, \qquad z \in \T,\: \alpha \in A.\]

\end{enumerate}

In this case, we write $\pt \sqsubset \Gamma \sqsubset \calU$.

We say $\Gamma$ is {\rm bounded above and below} or just {\rm bounded} in $\cbar$ with constant $K \ge 1$ if the mappings $\phi_\alpha$ can be chosen so that the resulting family 
$\Phi$ is bi-equicontinuous on $\T$ in the sense that we can find $K \ge 1$ such that 

\vspace{-.2cm}
\[ \frac{1}{K} \le |\phi_\alpha^\#(z)| \le K, \qquad z \in \T, \: \alpha \in A.\]

In this case, we write $\pt \sqsubset \Gamma \sqsubset \cbar$.

\end{definition}

\begin{proposition}[\cite{Com4} Proposition 6.1] Let $\calA = \{\Aaaa\}_{\alpha \in A}$ be a bounded family of pointed annuli, for each $\alpha$ let $\gamma_\alpha$ be the equator of $A_\alpha$ and let $\Gamma$ be the resulting family. Then $\pt \sqsubset \Gamma \sqsubset \calA$.
\end{proposition}

A result we shall need later is the following. 

\begin{theorem} [Transitivity of Boundedness, \cite{Com4} Theorem 6.1]$\begin{array}{cc} &\end{array}$\\
If $\pt \sqsubset \calU \sqsubset \calV \sqsubset \calW$ where $\calW$ is either bounded or all the domains of $\calW$ are $\cbar$, then 
$\pt \sqsubset \calU \sqsubset \calW$. This includes the degenerate case where $\calU$ is a family of simple closed curves.
\end{theorem}

We need to consider families of functions defined on bounded families of pointed discs or annuli which satisfy some suitable boundedness property which generalizes the notion of equicontinuity on compact subsets for a normal family of functions which are all defined on the same domain. We now give the precise definition we require.

\begin{definition}
Let $\calU = \{\Uaua\}_{\alpha \in A}$ be a non-degenerate family of hyperbolic pointed domains and  let $\calF = \{f_\alpha\}_{\alpha \in A}$ be a family of analytic functions 
with each $f_\alpha$ defined on $U_\alpha$. 

We say that 
$\calF$ is {\rm equicontinuous
on compact subsets of $\calU$} if for any hyperbolic distance $R >0$ there 
exists $M \ge 0$ depending on $R$ such that, for each $\alpha \in A$, $f_\alpha$ is $M$-Lipschitz with respect to the spherical metric within hyperbolic distance $\le R$ of $u_\alpha$ in $U_\alpha$. In this case we write $\calF \triangleleft \,\calU$.

We say that $\calF$ is {\rm bi-equicontinuous 
on compact subsets of $\calU$} if each $f_{\alpha}$ is a (locally injective) covering map onto its image and for each hyperbolic 
radius $R >0$ there exists $K(R) \ge 1$ such that, for each $\alpha \in A$, 
$f_\alpha$ is $K$-Lipschitz and locally $K$-bi-Lipschitz with respect to the spherical metric 
within hyperbolic distance $\le R$ of $u_\alpha$ in 
$U_\alpha$. In this case we write $\calF \bowtie \,\calU$.
\end{definition}

Note that if all the domains in our family are the same, then equicontinuity on compact subsets and local boundedness correspond with the standard 
definitions for a family of analytic functions on some domain (where we use the spherical rather than the Euclidean topology). In this 
context, bi-equicontinuity corresponds to having a normal family of 
covering maps where all of the normal limit functions must also be covering maps. In addition if all the functions in $\calF$ are univalent, then so are any limit functions.

We also need definitions of local boundedness, equicontinuity and bi-equicontinuity for the case of a family of functions defined on a bounded family of closed curves. 

\begin{definition}   
If $\Gamma$ is a bounded family of simple smooth closed curves as above and $\calF$ is a family of ${\mathrm C}^1$ functions defined on $\Gamma$, then we say that $\calF$ is  \emph{equicontinuous} on $\Gamma$ if there exists $K \ge 0$ such that for each $\alpha \in A$, $f_\alpha$ is $K$-Lipscitz with respect to the spherical metric on $\gamma_\alpha$. In this case we write $\calF \triangleleft \Gamma$. 

We say $\calF$ is \emph{bi-equicontinuous} on $\Gamma$ if there exists $K \ge 1$ such that for each $\alpha \in A$, $f_\alpha$ is a covering map from $\gamma_\alpha$ to its image $f_\alpha(\gamma_\alpha)$, $f_\alpha$ is $K$-Lipschitz and locally $K$-bi-Lipschitz on $\gamma_\alpha$. In this case we write $\calF \bowtie \Gamma$.

\end{definition}
 
For the sake of convenience we will use $\DO$ to refer to a constant family of pointed discs with base point $0$ where the index is suppressed. We 
remark in passing that the classical example of the family $\cal S$ of 
univalent functions on $\D$ which satisfy $f(0) = 0$, $f'(0) = 1$, is a good 
example of a bi-equicontinuous family on $\DO$.  However, it is important to note that we do not insist that all our functions be univalent and in practice many of the examples we will be working with will be covering mappings of bounded degree from one bounded family of annuli to another. 

Before we can go on to
state any more results, we need one further definition. 

\begin{definition} Let $\calU = \{\Uaua\}_{\alpha \in A}$ and $\calV = \{\Vava\}_{\alpha \in A}$ be families indexed by the same set $A$ with $\calU$ being non-degenerate. We say that a family $\calF = \{f_\alpha\}_{\alpha \in A}$ {\rm maps $\calU$ to $\calV$} if for each $\alpha$ $f_\alpha(U_\alpha) \subset V_\alpha$ and $f_\alpha(u_\alpha) = v_\alpha$ and we write $\calF: \calU \mapsto \calV$. 

By convention, if in addition each $f_\alpha$ is a covering map, we will require that $f_\alpha(U_\alpha) = V_\alpha$.
\end{definition}

\begin{theorem}{{\bf (Bi-Equicontinuity and Bounded Families of Discs\\ \cite{Com4} Theorem 4.1)}} Let $\calU = \{\Uaua\}_{\alpha \in A}$ be a family of pointed discs and let $\calF = \{f_\alpha\}_{\alpha \in A}$ be
the corresponding family of normalized inverse Riemann maps where for each 
$\alpha$, $f_\alpha$ maps $\D$ to $U_\alpha$ with $f_\alpha(0) = u_\alpha$ and $f_\alpha'(0) >0$. 
Then $\pt \sqsubset \calU \sqsubset \cbar$ if and only if 
$\calF$ is bi-equicontinuous on compact subsets of $\D$, i.e. $\calF \bowtie \DO$. 

\end{theorem}

We observe that the original theorem stated that for a general family $\calU$ of pointed hyperbolic domains where $\calF$ was then the family of suitably normalized universal covering maps, we had that $\calU$ is hyperbolically non-degenerate if and only if $\calF \bowtie \DO$. However, it is clear that in the case of families of simply connected domains, hyperbolic non-degeneracy is equivalent to boundedness (as the connectivity cannot decrease any further under \Car limits). Details about the more general version of this theorem can be found in \cite{Com4}.

The following is immediate both for bounded families of discs or annuli and for the case of families of simple closed pointed curves. 

\begin{lemma}[Composition of Equicontinuous Families] \hspace{1cm}\\
Suppose $\calU$ and $\calV$ are bounded, $\calF \triangleleft\, \calU$, $\calG \triangleleft \calV$, $\calF : \calU \mapsto \calV$ (with all families indexed by the same set) and let $\calG \circ \calF$ denote the family of composite maps (defined in the obvious way). Then $\calG \circ \calF \triangleleft \calU$. If in addition  
$\calF \bowtie \calU$ and $\calG \bowtie \calV$, then $\calG \circ \calF \bowtie \calU$.
\end{lemma}

We continue this section with the following result which shows
how bounded families of annuli may be obtained from 
bounded families of discs.

\begin{theorem}[\cite{Com4} Theorem 6.2]
Let $\calU$, $\calV$ be two bounded families of pointed discs with  
$\pt \sqsubset \calU \sqsubset \calV \sqsubset \cbar$. For each $\alpha$, let $A_\alpha$ be the conformal annulus 
$V_\alpha \backslash \overline {U_\alpha}$, let $\gamma_\alpha$ be the equator of $A_\alpha$ and let $a_\alpha$ be a point on $\gamma_\alpha$. Then $\calA = \{\Aaaa\}_{\alpha \in A}$ is a bounded family of pointed annuli.
\end{theorem}

Recognizing bounded families of annuli is important to us and so we include another result for this purpose. 

\begin{theorem}[\cite{Com4} Theorem 4.2] Let $\calA = \{(A_\alpha,a_\alpha)\}_{\alpha \in A}$ be a family of pointed annuli. Then $\calA$ is bounded if and only if we can find $\delta_1, \delta_2 >0$ and curves $\eta_\alpha$ for each $\alpha$ such that 

\begin{enumerate}

\item Each of the points $a_\alpha$ lies on $\eta_\alpha$ which separates the complementary components of $A_\alpha$ and 
the spherical distance from $\eta_\alpha$ to the complement of 
$A_\alpha$ is bounded below by $\delta_1$.

\item The spherical diameters of the complementary components of each $A_\alpha$ are $\ge \delta_2$.
\end{enumerate}
Furthermore, for any limit $\Aa$, the spherical distance to the boundary $\delta_A^\#(a)$ and the modulus $\bmod A$ are
bounded above and below, the bounds depending only on the constants $\delta_1$ and 
$\delta_2$.

\end{theorem}

\section{Polynomial-Like Mapping Sequences}

Having spent the last three sections developing tools, we can now finally 
turn to making use of these tools. We begin by recalling the classical 
definition of a polynomial-like mapping $\fuv$ of degree $d$ 
which consists of two simply connected open sets $U$ and $V$ with $U$ 
compactly contained in $V$, and a mapping $f$ which is a degree 
$d$ proper map of $U$ onto $V$. In our version of this we allow the 
mappings and also the sets they are defined on to vary at each stage of the 
iterative process and the definition is as follows. 

\begin{definition} Let $\calU = \{U_m, u_m\}_{m=0}^\infty$ and $\calV = \{V_m, u_m\}_{m=0}^\infty$ be two sequences of pointed discs with the same base points and for each $m \ge 0$, let $f_{m+1}$ be a mapping from $U_m$ onto $V_{m+1}$. We say that the sequence of triples $\fuvm$ is a 
{\rm polynomial-like mapping sequence} or simply {\rm polynomial-like sequence of degree bound $d$} if we can find a constant $K \ge 1$ such that the following are satisfied:

\begin{description}

\item[PL1] $\calU$ and 
$\calV$ are bounded sequences of pointed discs with $||\,\calU||, || \calV || \le K$  and $\pt \sqsubset \calU \sqsubset \calV$ (again with this constant $K$).
 
\item[PL2] For each $m \ge 0$, $f_{m+1}$ is a proper holomorphic map of degree\\ 
$2 \le d_{m+1} \le d$ of 
$U_m$ onto $V_{m+1}$ with $f_{m+1}'(u_m) = 0$.

\item[PL3] For each $m \ge 0$, the critical values of $f_{m+1}$ lie within hyperbolic\\ distance $\le K$ of $u_{m+1}$ in $V_{m+1}$.

\end{description}
\end{definition}

For convenience, we use the notation of Section 4 and write $\FUV$ to denote a polynomial-like sequence as above.

One simple consequence of the three properties above is that since by {\bf PL1}
$\calU$ is bounded, if $\phi_m$ is the inverse Riemann mapping from $\D$ to each $U_m$ with $\phi_m(0) = 0$, $\phi_m'(0) >0$, then by Theorem 4.3, if $\Phi$ denotes the corresponding family, then $\Phi \bowtie \Do$. From the boundedness of $\calV$ and the Schwarz lemma we must also have $\calF \circ \Phi \triangleleft \Do$. Combining these two observations, we have the following:

{\bf PL4} $\calF \triangleleft \,\calU$.  

Although this property is not logically independent of the other three, it is convenient for us to list it in this way.

If the sequences of sets $\{U_m\}_{m=0}^\infty, \{V_m\}_{m=0}^\infty,$ are constant as well as the 
sequence of mappings
$\{f_{m+1}\}_{m=0}^\infty,$, we recover the classical definition. It is worth making a remark at this point concerning
the set $V_0$. This set is not strictly necessary in the above definition. However, we will need it later in 
the proof of the non-autonomous versions of the straightening theorem and so we include this set as part
of the definition. 

Another more or less immediate consequence of this definition is a compactness property for pointwise limits. Suppose $\FUVn$ is a sequence each member of which is itself a polynomial-like sequence. We say $\FUVn$ converges \emph{pointwise} to another sequence $\FUV$ if for each fixed $m$, $(U_m^n, u_m^n) \to (U_m,u_m)$ and  $(V_m^n, u_m^n) \to (V_m,u_m)$ in the \Car topology while the functions $f_{m+1}$ converge uniformly on compact subsets of $U_m$ to $f_{m+1}$. The following result is an easy application of Cantor diagonalization and is left as a routine exercise for the reader. 

\begin{theorem} Let $\FUVa$ be a family of polynomial-like sequences as above and suppose that all members of this family have the same bounds $D$ and $K$. Then there exists a convergent subsequence which converges to a pointwise limit $\FUV$ which is another polynomial-like sequence with the same bound.
\end{theorem}

A similar result exists for sequences of polynomials (e.g. \cite{Com4}) and it is also related to the compactness result for polynomial-like mappings of fixed degree quoted in McMullen's book (\cite{McM} Page 72, Theorem 5.8).

As with sequences of polynomials, we will need 
some notation for the compositions of the functions involved and for such a sequence, we set
\vspace{.2cm}
\begin{eqnarray*}
F_m &=& f_m \circ \cdots \cdots \circ f_2 \circ f_1,\quad n \ge 1,\\  
F_{m,n} &=& f_n \circ \cdots \cdots \circ f_{m+2} \circ f_{m+1}, \quad 0 \le m \le n. 
\end{eqnarray*}

\vspace{-.2cm}
The classical concept of polynomial-like mappings
of course plays a central role in complex dynamics, especially in the areas of renormalization 
and the study of the structure of the Mandelbrot set. The starting point for these is the 
straightening theorem. 

\begin{theorem}(Douady, Hubbard) A polynomial-like map $\fuv$ is hybrid-equivalent to a polynomial. 
\end{theorem}

Recall that a hybrid equivalence is a quasiconformal equivalence whose dilatation is equal to zero on 
the filled Julia set associated with the polynomial-like mapping. In order to state the non-autonomous version
of this result, we will need to extend the classical terminology in a meaningful way to cover the 
non-autonomous case. We first need notions of equicontinuity and bi-equicontinuity which are appropriate for families of quasiconformal mappings. Again we let ${\mathrm d}^\#(x,y)$ denote the spherical distance between two points $x, y \in \cbar$. 

\begin{definition}
Let $\calU = \{\Uaua\}_{\alpha \in A}$ be a non-degenerate family of pointed discs or annuli and let $\Phi = \{\Phi_\alpha\}_{\alpha \in A}$ be a family of continuous functions defined on the domains of $\calU$. 

We say that $\Phi$ is {\rm equicontinuous} on $\calU$ if for every 
$\epsilon >0$, there exists $\delta >0$ such that for every $m \ge 0$, if $x$ and $y$ are points 
in $U_m$ with  
${\mathrm d}^\#(x,y)  < \delta$, then ${\mathrm d}^\#(\Phi_m(x), \Phi_m(y))< \epsilon$.  

We say that $\Phi$ is {\rm bi-equicontinuous} on $\calU$ if each $\Phi_\alpha$ is bijective, the image family $\{(\Phi_\alpha(U_\alpha), \Phi_\alpha(u_\alpha))\}_{\alpha \in A}$ is non-degenerate and the family $\Phi^{\circ -1}$ of inverse mappings is equicontinuous on this image family of pointed domains.

\end{definition}

\begin{definition} Let $\fuvm$ and $\fuvmt$ be two polynomial-like sequences with degree bound $d \ge 2$.
We say that two such sequences are 
{\rm quasiconformally conjugate} 
if there exists $K >1$, and a sequence $\Phi = \{\Phi_m\}_{m=0}^\infty$ of mappings  
for which we have the following:
\vspace{-.2cm}
\begin{enumerate}

\item For each $m \ge 0$, $\Phi_m$ is an orientation-preserving homeomorphism which maps 
$U_m$ onto $\tilde U_m$ and sends the base point 
$u_m$ of each $U_m$ to the corresponding 
base point $\tilde u_m$ of $\tilde U_m$;

\item For each $m \ge 0$, $\Phi_m$ (and hence $\Phi_m^{\circ -1}$) is $K$-quasiconformal;

\item The family  $\Phi$ is bi-equicontinuous 
on $\calU$. 

\item For each $m \ge 0$, $\Phi_{m+1} \circ f_{m+1} \circ \Phi_m^{\circ -1} = \tilde f_m$ on $\tilde U_m$. 

\end{enumerate}
\end{definition}

Note that in view of the boundeness of the families $\calU$, $\calUp$, Theorem 4.1 and Mori's Theorem (e.g. \cite{Ahl} page 47), we obtain locally  
uniform H\"older estimates on a quasiconformal equivalence as above. To be more precise, given $R >0$, there exist positive constants $C = C( R )$, $\alpha = \alpha( R )$ such that for every $m \ge 0$ and every $x, y \in U_m$ which are within hyperbolic distance $\le R$ of $u_m$, we have that  
$|\Phi_m(x) - \Phi_m(y)| < C|x - y|^\alpha$ with the family of inverse mappings also satisfying a 
similar condition.

Observer that this definition is considerably weaker than the classical one. However, 
in view of the bi-equicontinuity properties, it is strong enough to make meaningful comparisons 
between the dynamics of different non-autonomous dynamical systems. A similar notion has also been 
introduced by Kolyada and Snoha \cite{KS} in the context of real dynamical systems where it is referred
to as equiconjugacy. 

\begin{definition}
For  a polynomial-like sequence $\FUV$, we define the {\rm iterated filled Julia set} and {\rm iterated Julia set} at time $m$, $\Km, \Jm$, respectively by 
\[ \Km = \{z \in U_m: F_{m,n}(z) \notin V_n\setminus U_n,\:\:n > m\}, 
\qquad \Jm = \partial \Km .\]
\end{definition}

Note that it is trivial to show that the Julia sets 
$\Jm$ and the complements $U_m \setminus \Jm$ are completely invariant in the same sense as 
given in the statement of Theorem 1.1 for sequences of polynomials. 

With this behind us, we can now define our version of hybrid equivalence. 

\begin{definition}

We say that two polynomial-like sequences $\FUV$ and $\FUVt$ are 
{\rm hybrid-equivalent} if they are quasiconformally equivalent via a sequence $\Phi = \{\Phi_m\}_{m=0}^\infty$ where for each $m \ge 0$ the dilatation of each $\Phi_m$ 
is zero almost everywhere on $\Km$, the filled Julia set at time $m$ for 
$\FUV$. 
\end{definition}

Again, we remark that this is a weaker definition than the classical one, but still
strong enough to be useful. 

Given a classical  polynomial-like map $\fuv$, we may assume that the boundaries of $U$ and $V$
are smooth and in fact can be taken to be closed analytic arcs. The key to this useful trick is 
simply
to restrict the domain of $f$ to an appropriate set and it turns out that this idea will be 
very important to us in proving our version of the straightening theorem. This motivates the 
following definition. 

\begin{definition} If $\FUV = \fuvm$ is a polynomial-like sequence, we say the sequence
$\FUVp = \fuvmp$ is a {\rm polynomial-like restriction} (or simply {\rm restriction}) of $\FUV$ if, for each $m \ge 0$, 
$U'_m \subset U_m$, and $\FUVp$ is also   
polynomial-like and each $f_{m+1}$ has the same degree $d_{m+1}$ on $U_m'$ as on $U_m$.
\end{definition}

An important desirable feature of this definition is that making such a restriction does not change the iterated Julia sets. 

\begin{lemma} If $\FUVp$ is a restriction of the polynomial-like sequence $\FUV$, then the iterated filled Julia sets for $\FUVt$ are the same as those for $\FUV$. 
\end{lemma}

\proof Let $\Km$ and $\Kmp$ be the iterated filled Julia sets for the two sequences. Since the mappings for the two sequences have the same degrees, it follows that if $z \notin U'_i$ for some $i$, then $f_{i+1}(z) \notin V'_{i+1}$. Thus, for each $m$, 

\[ \Kmp = \bigcap_{n \ge m}{F_{m,n}^{\circ -1}(V'_n)} = 
\left (\bigcap_{n \ge m}{F_{m,n}^{\circ -1}(V_n)} \right )\cap U'_m = \Km \cap U'_m.\]
 
 Suppose for now that $m \ge 1$. Then there can be no points of $\Km$ in $V'_m \setminus U'_m$ since the preimages of a point in $\K_m \cap 
(V'_m \setminus U'_m)$ would all lie in $U'_{m-1}$ and thus in 
$\K'_{m-1}$ while $f_m$ maps $\K'_{m-1}$ to $\Kmp \subset U'_m$. 
Suppose now that we have a point which lies in $\Km \setminus \Kmp$. 
Thus $\Km \subset U'_m \cup (U_m \setminus V'_m)$, and since $U'_m \cup (U_m \setminus V'_m)$ is disconnected, it follows by the 
definition of $\Km$, that we can find a smallest $n$ such that  $F_{m,n}^{\circ -1} (V_n)$ has at least one component which meets $U_m \setminus V'_m$ but not $U'_m$. However, since $V_n \supset V'_n$ and $f_{i+1}$ has the same degree on $U'_i$ as on $U_i$ for each $m \le i <n$, it follows by an easy induction that every component of the inverse image of $F_{m,n}^{\circ -1} (V_n)$ must meet $U'_m$. With this contradiction, we obtain that $\Kmp = \Km$ for each $m \ge 1$. Finally the case $m = 0$ follows by taking inverse images under $f_1$ and remembering that this mapping has the same degree on $U_0'$ as on $U_0$.  $\Box$

We now state the non-autonomous version of the straightening theorem. 

\begin{theorem}
Let $\FUV$ be a polynomial-like mapping sequence with degree bound $d$. Then there exists a 
polynomial-like restriction $\FUVp$ and a bounded monic centered sequence of polynomials $\Pm$ such 
that $\FUVp$ is hybrid-equivalent to $\Pm$. Moreover, the degree of each $P_{m+1}$ is $d_{m+1}$, the degree of
$f_{m+1}$.
\end{theorem}

%%%%%%%%%%%%%%%

\section{Construction of the Mating}

We are now in a position to commence 
proving Theorem 5.3. The proof will
consist of three main parts. The first part is the construction of the desired polynomial-like restriction. The second part is a uniformization whose point is to ensure that we will be dealing
with  a family of functions defined on essentially the same set where it will be easy to see that we obtain a bi-equicontinuous family of quasiconformal maps of uniformly bounded 
dilatation. The final part is the construction of the required conjugacy which will be done by a non-autonomous version of the usual fundamental domains argument. We 
now begin with the first stage of Part 1. 

\vspace{.4cm}
{\bf Part 1: Construction of the Polynomial-like Restriction}

\vspace{.2cm}
{\bf Stage I:} Restriction of the Domains. 

We begin by noting that by Theorem 4.1 we can perform a (non-autonomous) conjugacy if necessary which 
consists of a uniformly bi-Lipschitz sequence of M\"obius transformations so that we may assume that 
$u_m = 0$ and $\infty \notin V_m$ for each $m \ge 0$ where we recall that have assumed that 
$u_m$ is a critical point of $f_{m+1}$.
Let $K$ be the number appearing in Definition 5.1 for polynomial-like sequences.  Let $B > 1$ and let $T_m$ be the conformal annulus 

\[T_m = \{z \in V_m: BK < \rho_{V_m}(z\,,0) < 2BK\} \]

where $\rho_{V_m}(\cdot \,,\cdot)$ denotes hyperbolic distance in $V_m$.
Let $\Gm$ be the equator of this annulus and let $V_m'$ be the topological disc enclosed by this curve. Let $\calT$ and $\calVp$ 
be the families of pointed annuli and pointed discs obtained in this way where we choose the base point $t_m$ of each annulus $T_m$ to be an arbitrary point on the equator while all the discs have base point $0$. Finally, let $\bfG$ denote the degenerate family of pointed curves obtained from the curves $\Gamma_m$ with the base points $t_m$. Note that, since $\infty$ does not lie in any of the domains $V_m$, by Corollary 2.1 all the sets $V_m'$, $T_m$  lie within uniformly bounded Euclidean distance of $0$.  From this it follows easily that the estimates on the Lipschitz constants for all the statements concerning equicontinuity and bi-equicontinuity which we will need in what follows can be obtained by estimating ordinary derivatives. 

Now for each $m \ge 0$ let $S_m$ be the inverse image under $f_{m+1}$ of  $T_{m+1}$. Since all the critical values of $f_{m+1}$ lie in the bounded complementary component of $T_{m+1}$, $S_m$ will also be a conformal annulus and we denote the resultant family by $\calS$ where the base point $s_m$ of each $S_m$ is an arbitrarily chosen preimage of $t_{m+1}$ under $f_{m+1}$. Let $\gm$ be the equator of $S_m$ (note that $s_m$ will lie on $\gm$) and let $U_m'$ be the topological disc enclosed by $\gm$. Finally, let $\calUp$ be the resultant family of pointed discs where each disc again has base point $0$, (which is still a critical point of each $f_{m+1}$) and let $\bfg$ denote the degenerate family obtained from the curves $\gamma_m$ with base points $s_m$. 

Denote the new sequence $\fuvmp$ by $\FUVp$.

\begin{description}

\item[I.1] By choosing $B$ suitably and making $K$ larger if needed, 
$\FUVp$ is a polynomial-like restriction of $\FUV$.

First we note that since $\calU \sqsubset \calV$, we can choose $B$ large enough so that $\calU \sqsubset \calV' \sqsubset \calV$. By Lemma 2.1, 
$\pt \sqsubset \calU' \sqsubset \calU$ whence by Theorem 4.2 on the transitivity of boundedness we have
$\pt \sqsubset \calU' \sqsubset \calV'$, thus proving {\bf PL1} for $\FUVp$ (note that we can choose $B$ large enough so that $||\,\calU'||, ||\calV'|| < 2K$).

{\bf PL2} Follows simply from the fact that $\Gamma_{m+1}$ contains all the critical values of $f_{m+1}$ in view of our choice of $B$. 

Recall the standard estimates estimates given in Section 2 for the hyperbolic metric in terms of Euclidean distance to the boundary for simply connected domains. Using this and {\bf PL3} for the original sequence $\FUV$ combined with our choice of the curves $\Gm$ 
tells us that we can connect the critical values of each $f_{m+1}$ to $u_{m+1}$ using curves which are uniformly bounded away (in terms of the Euclidean distance between sets) from the curves $\Gamma_{m+1}$. Using the same estimates again, we obtain {\bf PL3} for $\FUVp$ (as with {\bf PL1}, we can make $B$ large enough so that this condition is satisfied with bound $2K$).

\item[I.2] $\pt \sqsubset \calS \sqsubset \cbar$ and 
$\pt \sqsubset \calT \sqsubset \cbar$.

Dealing first with $\calT$, since $\pt \sqsubset \calV \sqsubset \cbar$, the diameters of the unbounded complementary components are bounded away from $0$.  Since $\calUp$ is bounded, if we apply Theorem 4.3 and the Koebe one-quarter theorem to the family of suitably normalized inverse Riemann mappings from $\D$ to the discs $U_m'$, we can say the same about the bounded complementary components. On the other hand, Koebe and the distortion estimates for univalent mappings (e.g. \cite{CG} page 3, Theorem 1.6) applied to the boundedness of $\calV$ also tell us that 
the spherical (and hence Euclidean) distance of the curves $\Gm$ to the boundaries $\partial T_m$ is uniformly bounded away from $0$. The boundedness of $\calT$ follows by Theorem 4.5. 

Turning now to $\calS$, since $\pt \sqsubset \calU \sqsubset \cbar$, again the diameters of the unbounded complementary components are clearly bounded away from $0$. If we now apply the Schwarz lemma for the hyperbolic metric to the mappings $f_{m+1} : U_m \mapsto V_{m+1}$, and apply Koebe and Theorem 4.3 as before, we can also conclude that the diameters of the bounded complementary components will be bounded away from $0$. By Lemma 2.1, the equator $\gm$ of $S_m$ is contained within uniformly bounded hyperbolic distance of $0$ in $U_m$ while, by the Schwarz lemma, its hyperbolic distance in $U_m$ from $\partial S_m$ is bounded below. The same argument using the distortion theorems as in the last paragraph then shows that 
the spherical and Euclidean distances from the curves $\gm$ to the boundaries $\partial S_m$ are uniformly bounded 
away from $0$. Hence $\calS$ is also bounded, again by Theorem 4.5. 

\item[I.3] $\pt \sqsubset {\boldsymbol \gamma} \sqsubset \calS$, 
$\pt \sqsubset \boldsymbol \Gamma \sqsubset \calT$ and $\pt \sqsubset {\boldsymbol \gamma} \sqsubset \cbar$, $\pt \sqsubset \boldsymbol \Gamma \sqsubset \cbar$.

The first pair of statements follow from applying Proposition 4.1 to {\bf I.2} above. 
The second pair then follow from the boundedness of $\calS$ and $\calT$ above and 
Theorem 4.2 on the transitivity of boundedness. 

\item[I.4] $\calF \bowtie \calS$. 

We saw in the proof of {\bf I.2} that the equators of the annuli of the families $\calS$ and $\calT$ are uniformly bounded away from the boundaries of these annuli and applying Koebe then gives us the desired estimates on the derivatives and hence the local (spherical) Lipschitz constants at the base points. By Theorem 4.1 the moduli of these annuli are uniformly bounded above and so by Corollary 2.1, the constants $c$ in the lower bound on the hyperbolic metric for these annuli are uniformly bounded below away from $0$. Hence, given any hyperbolic radius $r$, we can find $\delta ( r ) > 0$ such that 
for all $m$, if $z \in T_m$ is of hyperbolic distance $\le r$ from the base point $t_m$ of this annulus, then $\delta (z) \ge \delta( r )$. If we then connect $z$ to $t_m$ along a hyperbolic geodesic using a suitable chain of discs of radius $\tfrac{\delta( r )}{2}$ (so that the discs with the same centres and twice the radii still lie in $T_m$ and the number of such discs we require is uniformly bounded with respect to $r$) and apply the distortion theorems for univalent mappings on each such disc in turn, we obtain the desired estimates on the Lipschitz constants everywhere else on these annuli. 

\end{description}

\vspace{.2cm}
{\bf Part 2: Uniformization of the Domains for $\FUVp$}

\vspace{.2cm}
{\bf Stage II:} Uniformization of the Outer Domains for $\calVp$

For each $m \ge 0$ let $W_m$ be the topological disc $\cbar \setminus \overline{V'_m}$ and let $\calW$ be the family resulting from the pointed discs $(W_m, \infty)$ (recall that we have assumed that $\infty \notin V_m$). Let $R >1$ and for each $m$ let $\psi_m$ be the Riemann mapping 
from $W_m$ to $\cbar \setminus \overline{D(0,R)}$ with $\psi_m'(\infty) > 0$. Let $\Psi$ be the family of functions we obtain in this way.

{\bf Remarks}

\begin{description}
\item[II.1] $\Psi \bowtie \calW$.

The mappings $\psi_m$ send $\infty$ to itself and, since $\calVp \sqsubset \calV \sqsubset \cbar$ (or, as we remarked earlier, the discs $V_m'$ are contained within uniformly bounded Euclidean distance of $0$), the absolute values of the derivatives $\psi_m'(\infty)$ are uniformly bounded above and below away from $0$ by Koebe. The result now follows from the distortion theorems for univalent mappings.

\item[II.2] We can extend the domains of the mappings $\psi_m$ univalently so that $\Psi \bowtie \calT$. 

The extension to the domains of $\calT$ follows immediately by the results on Schwarz reflection across analytic arcs (e.g. \cite{Lang} Chapter IX, Theorem 2.2) and it is easy to see that these extensions remain univalent. Since, as we saw in the course of proving {\bf I.2}, the equators $\Gamma_m$ of each $T_m$ are uniformly bounded away from $\partial T_m$, the estimates on the hyperbolic metric for simply connected domains in Section 2, show that that for each of the composite domains $(\cbar \setminus \overline {V'_m}) \cup T_m$, the hyperbolic distance from $\infty$ to points on $\Gamma_m$ is uniformly bounded above.  The conclusion then follows from the distortion theorems by a similar argument to that for estimating the Lipschitz constants for $\calF$ on $\calS$ in the proof of {\bf I.4}.
\end{description}

{\bf Stage III:} Mapping to Round Annuli

Let $\cpsi_m$ be a Riemann map which maps the annulus 
$L_m := V'_m \setminus \overline {U'_m}$ to the round concentric annulus $\cLm = A(0,R_m,R)$ where $R$ is as above and $R_m$ is such that this round annulus has the same modulus as $L_m$. Let us denote by $\cgm$ and $\cGm$ the inner and outer boundary curves of $\cLm$. 

Let $l_m$ and $\clm$ be base points for each annulus $L_m$ and $\cLm$ respectively which are on the equators of these annuli with $\cpsi_m (l_m) = \clm$. Finally, let $\calL$ and $\calLc$ be the two families we obtain with these base points and let $\Psic$ denote the resulting family of Riemann maps defined on $\calL$.

{\bf Remarks}
\begin{description}

\item[III.1] $\pt \sqsubset \calL \sqsubset \cbar$ and $\pt \sqsubset \calLc \sqsubset \cbar$.

The boundedness of $\calL$ is a consequence of {\bf PL1} for $\FUVp$ and Theorem 4.4. The boundedness of $\calLc$ then follows easily from this and Theorem 4.1  
as the annuli of this family are round and concentric and all have the same outer radius while their moduli are uniformly bounded above and below in view of the boundedness of $\calL$ and conformal invariance.

\item[III.2] $\Psic \bowtie \calL$.

The annuli of $\calLc$ are round, have the same outer radius and from above have moduli which are uniformly bounded above and below. The estimates on the derivatives and hence the local Lipschitz constants at the base points follow as before from the Koebe one-quarter theorem while Corollary 2.1 and the distortion theorems give us the estimates elsewhere as in the proofs of {\bf I.4} and {\bf II.2} (this conclusion can also be obtained more quickly by appealing directly to \cite{Com4} Theorem 4.2).  

\end{description}

Now let $M_m = S_m \cup L_m \cup T_m$. Then $M_m$ is clearly a conformal annulus and we let $\calM = \{(M_m, t_m)\}_{m=0}^\infty$ where $t_m \in \Gm$ was the base point we chose for the pointed annulus $(T_m,t_m)$ introduced earlier. 

Again, we can use Schwarz reflection, this time to univalently extend the domain of definition of the functions $\cpsi_m$ to all of $M_m$. 

\begin{description}

\item[III.3] $\pt \sqsubset \calM \sqsubset \cbar$ and $\Psic \bowtie \calM$. Hence $\Psic \bowtie \calS$ and $\Psic \bowtie \calT$.

The inner complementary component of $M_m$ is that of $S_m$ while the outer complementary component of $M_m$ is that of $T_m$. By {\bf I.2} $\calS$ and $\calT$ are bounded families and we already saw that the (spherical) diameters of these components are bounded below. On the other hand, $M_m$ contains $T_m$ and its equator $\Gm$ and so, as we saw in proving {\bf I.2}, these complementary components are bounded away from $\Gm$ (on which our base point $t_m$ is located). The boundedness of $\calM$ now follows from Theorem 4.5.  

The bi-equicontinuity of $\Psic$ on $\calM$ 
follows from {\bf III.2} and a similar argument to that for $\Psi$ on $\calT$ in the proof of {\bf II.2} (note that the bound on the annuli of $\calM$ follows from the boundedness of this family above and Theorem 4.1). The bi-equicontinuity of $\Psic$ on $\calT$ is then immediate using the Schwarz lemma for the hyperbolic metric. On the other hand, the bound on the moduli of the annuli of $\calM$ and the estimates of Corollary 2.1 on the hyperbolic metric tell us that the curves $\gamma_m$ are uniformly bounded distance from the curves $\Gamma_m$ in the hyperbolic metric of $M_m$. The local boundedness and  bi-equicontinuity of $\Psic$ on $\calS$ then follow easily from this and the corresponding properties for $\calM$ (and once more the Schwarz lemma for the hyperbolic metric).

\end{description}

\vspace{.2cm}
{\bf Part 3: Construction of the Hybrid Conjugacy}

\vspace{.2cm}
{\bf Stage IV:} Extending the mappings $\psi_m$ quasiconformally to $\cbar \setminus U'_m$.

We want to extend the mappings $\psi_m$ to all of the sets $\cbar \setminus U'_m$ in a way that combines the dynamics of the functions
$f_{m+1}$ with that of  the mappings $z^{d_{m+1}}$. We have no hope of doing this conformally. However, it can be done in a quasiconformal fashion so as to give us a quasi-regular sequence. In order to do this, we will need to redefine $\psi_m$ on $V'_m \cap T_m$ where it has already been defined. We will be able to do this in such a way that the extended function will still be continuous. 

We begin by defining a sequence of functions $\check f_m$ which give rise to a family $\bfcF$ on the family $\bfcg$ arising from the curves $\cgm = C(0, R_m)$ (note that since the boundaries of the annuli $L_m$ are smooth the Riemann mappings $\cpsi_m$ extend continuously to the boundary and so we can take our base points to the be the images of the points $s_m$ under these mappings). Similarly we can define the family $\bfcG$ arising from the curves $\check \Gamma_m = C(0,R)$ whose base points are the images of the points $t_m$ under the mappings $\psi_m$ in the same way as above. 

The following is immediate from {\bf III.1} and Proposition 4.1 and Theorem 4.2.

{\bf IV.1}  $\pt \sqsubset \bfcg \sqsubset \cbar$ and $\pt \sqsubset \boldsymbol \bfcG  \sqsubset \cbar$.

\vspace{.2cm}
We next set

\vspace{-.2cm}
\[ \check f_m = \cpsi_{m+1} \circ f_{m+1} \circ \cpsi_m^{\circ -1}.\]
  
Since $\cpsi_m$ maps $\gamma_m$ to $\tgm$ and $\Gamma_m$ to $\tGm$, by {\bf I.4} and {\bf III.3} we obtain the following.

\begin{description}

\item[IV.2]  $\bfcF$ is a bi-equicontinuous family of ${\mathrm C}^1$ mappings from 
$\bfcg$ to $\bfcG$.

\end{description}

Let $\tRm = R^{1/d_{m+1}}$ and let $\tLm$ be the annulus $A(0,\tRm,R)$ and let $\calLt$ be the resulting family (note that the base points are not yet fixed, but we will be able to do this fairly soon).

Now let $\tgm$ and $\tGm$ be the inner and outer boundary curves of $\tLm$, and let $\bftg$ and $\bftG$ be the resulting families where again we define the base points later. 

The next stage is to construct suitable mappings $\ccpsi_m$ from $\cLm$ to $\tLm$ so that setting $\psi_m = \ccpsi_m \circ \cpsi_m$ will yield the desired extension of $\psi_m$. 

To begin with we define $\ccpsi_m$ on $\cGm$ by 

\vspace{-.2cm}
\[\ccpsi_m(z) = \psi_m \circ \cpsi_m^{\circ -1}(z)\]

and letting the base points for $\bfcG$ to be the images of the points $t_m$ under the mappings $\psi_m$. By {\bf II.2} and {\bf III.3} we immediately have that 

\begin{description}

\item[IV.3]  $\bfccPsi$ is a bi-equicontinuous family of sense-preserving bijective ${\mathrm C}^1$\\ mappings from $\bfcG$ to $\bftG$.

\end{description}

We now define the mappings $\ccpsi_m$ on the sets $\tgm$ by the following lifting diagram which also allows us to specify the base points for the family $\bftg$.

\vspace{-.6cm}
\begin{diagram} 
\cgm &&\rTo^{\ccpsi_m} && \tgm\\\\
\dTo^{\check f_{m+1}} & & & & \dTo_{z^{d_{m+1}}}\\ \\
\check \Gamma_{m+1} &&\rTo^{{\check{\check \psi}}_{m+1}} &&{\tilde \Gamma}_{m+1}
\end{diagram}

\vspace{.2cm}
Using {\bf IV.2} and {\bf IV.3}, we see that 

\begin{description}

\item[IV.4]  $\Psicc$ is a bi-equicontinuous family of sense-preserving bijective ${\mathrm C}^1$\\ mappings from $\bfcg$ to $\bftg$.

\end{description}

In \cite{Leh}, Lehto proves a Lemma in the course of which he considers the interpolation problem of finding a quasiconformal mapping from a round concentric annulus to itself which extends a given ${\mathrm C}^1$ identification of their boundaries. In our case, the annuli $\tLm$, $\cLm$ need not have the same moduli. However, by {\bf III.1} and Theorem 4.1, the moduli of the annuli $\cLm$ are bounded above and below. Since the degrees $d_{m+1}$ are uniformly bounded, the same is obviously true of the annuli $\tLm$. Then by first applying a suitable sequence of bi-equicontinuous quasiconformal mappings of uniformly bounded dilatation (e.g. an affine stretching in logarithmic coordinates), we can then assume without loss of generality that we are in a position to make use of Lehto's result. Applying {\bf IV.3} and {\bf IV.4} to the estimates on the partial derivatives and the complex dilatation in Lehto's proof, we obtain the following. 

\begin{description}

\item[IV.5] By making $K$ larger if necessary, the mapping $\ccpsi_m$ can be extended to a homeomorphism of $\overline {\check L}_m$ to $\overline {\tilde L}_m$ to give a bi-equicontinuous family of $K$-quasiconformal mappings.

\end{description}

Now define $\psi_m$ on all of $\cbar \setminus U'_m$ by

\vspace{-.1cm}
\[\psi_m(z) =  \left \{ \begin{array}{r@{\quad , \quad}l}
\ccpsi_m \circ \cpsi_m (z) & z \in {\overline V'}_m \setminus U'_m\\
\psi_m (z) & z \in \cbar \setminus {\overline V'}_m
\end{array} \right .\]
\vspace{.1cm}

Note that the two definitions agree on $\partial V_m$. Recall that we were originally able to 
extend the mappings $\psi_m$ conformally beyond the sets $\partial V_m$ using Schwarz reflection ({\bf III.3}). However, we needed to redefine these mappings on $V'_m \cap T_m$ to get our desired quasiconformal extension. If we then compare these two definitions and use Rickman's Lemma (e.g. \cite{DH} page 303 Lemma 2) on suitable neighbourhoods of the compact sets $\overline T_m \setminus V'_m$, we see that the above extensions of the functions $\psi_m$ to the sets $\cbar \setminus U'_m$ will be $K$-quasiconformal. 
By {\bf II.1}, {\bf III.3} and {\bf IV.5}, if we still let $\Psi$ denote the correpsonding family, we get the following. 

\begin{description}

\item[IV.6] $\Psi$ is a bi-equicontinuous family of $K$-quasiconformal sense-preserving homeomorphisms on the family $\{(\cbar \setminus U'_m, \infty)\}_{m=0}^\infty$ (not just on compact subsets).

\end{description}

{\bf Part V} Extension to a quasi-regular sequence of mappings on $\cbar$. 

For each $m \ge 0$ define a mapping $g_m(z)$ from $\cbar$ to itself by 
\vspace{.1cm}
\[ g_{m+1}(z) = \left \{ \begin{array}{r@{\: , \quad}l}
 f_{m+1}(z) & z \in U_m\\
\psi_{m+1}^{\circ -1} \left ( ( \psi_m (z))^{d_{m+1}} \right ) & z \in \cbar \backslash U_m
\end{array} \right . .\]
\vspace{-.1cm}

Note that by the definitions of $\cpsi_m$, $\ccpsi_m$, the mapping $g_{m+1}$ will be continuous on $\partial U'_m$ and thus everywhere on $\cbar$ for each $m \ge 0$.
For $0 \le m \le n$, let $H_{m,n} = 
g_n \circ \cdots \circ g_{m+1}$ be the compositions similar to those for polynomial-like mapping sequences introduced at the start of Section 5. Then for any given $m$, the interiors of the sets $H_{m,n}^{\circ -1}(L_n), \, n \ge m$ are pairwise disjoint. Another way of saying this is that
the sets $L_m \setminus \gm$ give us a set of fundamental domains for the non-autonomous dynamical system arising from the functions $g_{m+1}$ defined on the constant sequence of sets $\cbar$. Combining this observation with 
{\bf IV.5} above, we get the following. 

\begin{description}
\item[V.1] {\it The sequence $\{g_m\}_{m=1}^\infty$ is $K$-quasi-regular.} 
\end{description}

{\bf Part VI} Construction of the invariant ellipse field and the hybrid conjugacy. 

This part of the argument is relatively standard. Start by recalling that for 
a given function $f$,  $\mu_f$ is the standard notation for the complex 
dilatation of $f$. For each $n \ge 0$, define a Beltrami coefficient $\mu_m(z)$ 
(which corresponds to a change of complex structure) on $\cbar$ by 

\vspace{-.1cm}
\[\mu_m(z) = \left \{ \begin{array}{r@{\: , \quad}l} \mu_{H_{m,n +1}}(z) & 
H_{m,n}(z) \in L_n \setminus \gamma_n, \quad \mbox{for some} \quad n \ge m\\ 
0 & \mbox{otherwise}
\end{array} \right .  .\]

By the quasi-regularity of the sequence of mappings 
$\{\tilde g_m\}_{m=1}^\infty$, it follows
easily that $||\mu_m||_\infty < (K-1)/(K+1)$. Hence for each $m$ there exists a 
unique normalized solution $\Phi_m$ which has dilatation $\mu_m$ (almost everywhere) and fixes $0$, $1$ and $\infty$. It is not hard to check that for each $m \ge 0$, 
(e.g. using the calculations in \cite{Ahl} P. 9) that the mappings 
$P_{m+1} = \Phi_{m+1} \circ g_{m+1} \circ {\Phi_m}^{\circ -1}$ are analytic on all 
of $\cbar$. To be more precise, one can verify that the above definition of the Beltrami coefficients $\mu_m$ together with the fact that the mappings $g_{m+1}$ are analytic except (possibly) on the sets $\overline {L_m}$ give sufficient conditions for the mappings $P_{m+1}$ to be analytic. Finally, it is easy to see that these functions must have topological degree $d_{m+1}$
with a pole of order $d_{m+1}$ at $\infty$ and no other poles. Hence they must be polynomials. 

Now for each $m$ the support of $\mu_m$ is uniformly bounded in $m$. Hence
$\Phi_m$ is a K-quasiconformal orientation-preserving homeomorphism of $\cbar$ 
which by the equations (5) and (7) on page 92 of \cite{Ahl} has the form  
$\Phi_m(z) = z + {\cal O}(1/z)$ near $\infty$.  Combining this estimate at $\infty$ with the fact that the family $\Psi$ is bi-equicontinuous on $\calW$ by {\bf II.1}, we see that the leading coefficients of the polynomials are uniformly bounded above and below in absolute value (and must in fact be positive).

We now turn to estimating the other coefficients of our sequence. The Beltrami differential of each $\Phi_m$ is bounded above by $K$ and the supports of these differentials are uniformly bounded.  Also, the same is clearly true for the family of inverse maps $\Phi_m^ {\circ -1}$. By Theorem 7.3 on Page 21 of \cite{CG} page 21, the family arising from the mappings $\Phi_m$ is then bi-equicontinuous on $\cbar$. Since $\Phi_m(\infty) = \infty$ for each $m \ge 0$, combining this with the bi-equicontinuity of the family $\Psi$ on $\{(\cbar \setminus U'_m, \infty)\}_{m=0}^\infty$, we see that the sequence $\{P_m\}_{m=1}^\infty$ is uniformly bounded on any fixed large disc centered about $0$. By Cauchy's estimates, all the lower order coefficients will then be uniformly bounded as desired.

To complete the proof, we need to show that we have a 
hybrid conjugacy between $\fuvmp$ and $\Pm$. To do this we need to check that the mappings $\Phi_m$ satisfy properties \emph{1.}-\emph{4.} in definition 5.3 for quasiconformal equivalence and that the dilatations of these mappings are zero almost everywhere on the filled Julia sets. First note that we have already seen above that the functions $\Phi_m$ give rise to a bi-equicontinuous family in the sense given in Definition 5.2 which then satisfies property {\emph 3.} of Definition 5.3. Also, the dilatation $\mu_m$ of the mappings $\Phi_m$ is automatically zero on the filled Julia sets $\tilde K_m$ for the sequence
$\fuvmp$. The remaining properties follow immediately and with this the proof is finally complete. $\Box$

%%%%%%%%%%%%%%%

\section{Applications}

We present two applications of the results we have proven; one a generalization of a classical theorem, the other a counterexample. 

Recall the classical fact that, for a quadratic polynomial $P(z) = z^2 + c$, the Julia set is a quasicircle if $|c| < 1/4$. Our theory of polynomial-like sequences allows us to very quickly prove a non-autonomous version of this result.

\begin{theorem}
For a sequence $\{z^2 + c_m\}_{m=1}^\infty$ of quadratic polynomials with $||\{c_m\}_{m=1}^\infty||_{l^\infty} < 1/4$, the iterated Julia sets are quasicircles. 
\end{theorem}

\vspace{.2cm}
\proof For each $m \ge 0$, let $V_m = \cbar \setminus \overline {\mathrm D}(c_m, \tfrac{1}{4})$, $U_m =  P_{m+1}^{\circ -1}(V_{m+1})$ and make these into two families of pointed discs by letting all the base points be $\infty$. It is then easy to check that this gives rise to a polynomial-like sequence. The result then follows directly from the straightening theorem (Theorem 5.3).  $\Box$

This result was already proved by Br\"uck \cite{Br1} who obtained it by showing that the inside and outside of the Julia set were both John domains. There is essentially no loss of generality in considering sequences of this type as it was shown in \cite{Com1} that any bounded sequence of quadratic polynomials is conformally conjugate to a sequence of monic centered quadratic polynomials.

Recall the classical result that all quadratic-like mappings with disconnected Julia sets are hybrid equivalent. In the non-autonomous case this is no longer true. 

\begin{theorem}
There exists a bounded sequence of quadratic polynomials all of whose critical points escape and all of whose iterated Julia sets are totally disconnected, but which is not hybrid equivalent to any quadratic polynomial with disconnected Julia set (viewed as a constant sequence).
\end{theorem}

\proof 

For each $n \ge 1$ let $\Pmn$ be the following sequence defined by 
\[ P_m^n = \left \{ \begin{array}{r@{\quad , \quad}l}
(z-3)^2 & \quad m = (j+1)(j+2)/2 - 1 \mbox{\quad for some \quad} 1 \le j \le n\\
z^2  & \mbox{otherwise}  \end{array} \right.\]
In other words, we iterate once with $z^2$ followed by $(z-3)^2$, then twice
with $z^2$ followed by $(z-3)^2$, then three times with $z^2$ followed by
$(z-3)^2$ and then four times with $z^2$ and so on until $n$ of these steps have been 
carried out after which the sequence is simply $z^2$. 

Let $\Jmn$ denote the iterated Julia sets for these sequences. 
It is fairly easy to check that at the times 
$(j+1)(j+2)/2 - 2$, which are those times just before we iterate with $(z-3)^2$, $\J_{(j+1)(j+2)/2 - 2}^n \subset \cbar \setminus \overline D(3,1)$ 
while at all other times $\Jmn \subset \cbar \setminus \overline \D$. 
From this it follows easily that $|{P_m^n}'| > 2$ on $\Jmn$ for each $m \ge 1$.

The sequences $\Pmn$ are thus uniformly hyperbolic and also uniformly bounded in terms of degrees and coefficients and converge to the pointwise limit $\Pm$ where 

\vspace{-.4cm}
\[ P_m = \left \{ \begin{array}{r@{\quad , \quad}l}
(z-3)^2 & \quad m = (j+1)(j+2)/2 - 1 \mbox{\quad for some \quad} j \ge 1\\
z^2  & \mbox{otherwise}  \end{array} \right.\]

Appealing to the theory of hyperbolic non-autonomous Julia sets as developed in \cite{Com2}, we see from Theorem 3.4 in that paper that $\Pm$ is also hyperbolic and that the iterated Julia sets $\Jmn$ for the sequences $\Pmn$ converge in the Hausdorff topology converge to the corresponding Julia sets $\Jm$ for $\Pm$ (this also follows from a result of Sumi (\cite{Sum} page 583 Theorem 2.14) as well as a result in the paper of Sester (\cite{Ses} page 411, Proposition 4.1), both of whom were working in the context of polynomials fibered over a compact set). It is also easy to see that all the critical points escape to infinity and since the sequence is hyperbolic there can be no bounded Fatou components in view of Lemma 3.2 from \cite{Com2}. 

In fact, if we look at the times $(j+1)(j+2)/2 - 1$ which are precisely those times just before we start the cycles of iteration with $z^2$, the Julia set consists of $2^{j+1}$ congruent pieces of diameter $\sim 2^{-j-1}$ spaced around a circle about $0$ whose radius tends to $1$ as $j$ tends to infinity. Also, the separation between any two adjacent such pieces is also $\sim 2^{-j-1}$.

Now consider a standard example of a quadratic polynomial with disconnected Julia set e.g. $P(z) = z^2 + 2$. It is easy to check that, for this polynomial, the points on the circle ${\mathrm C}(0,2)$ escape to infinity and so there is a disc of radius $R < 2$ about $0$ which contains the Julia set. The preimages of this disc give two sets which are $\delta: = 2\sqrt {2-R} > 0$ apart and so the 
Julia set also consists of two pieces each of which is mapped onto the whole Julia set by $P$ and which are at least $\delta$ apart. 
It then follows that any supposed quasiconformal conjugacy between $\Pm$ and $\{P\}_{m=1}^\infty$ will map two points of distance $\sim 2^{-j}$ apart to two points which are distance $\ge \delta$ apart. As $j \to \infty$, this clearly violates the bi-equicontinuity requirement for any conjugacy and with this the result follows. $\Box$

\end{document}